\documentclass[11pt]{article}
\def\R{\ifmmode{\rm I\mkern-3.1mu
R\mkern1mu}\else{\rm I\kern-.18em  R\hskip1pt\
}\fi\relax}

\def\b{\beta}

\def\G{\Gamma}

\def\t{\tau}

\def\th{\theta}
\def\l{\lambda}
\def\L{\Lambda}

\def\s{\sigma}
\def\sou{\overline}
\def\so{\underline}

\def\f{\rightarrow}
\def\q{\forall}

\def\v{\vdash}

\def\mats{\ifmmode{ {\hbox{\bigreek s}} }\else{ {\bigreek s} }\fi\relax}
\def\matsin{\ifmmode{ {\hbox{\smgreek s}} }\else{
{\smgreek s} }\fi\relax}
\def\matt{\ifmmode{ {\hbox{\bigreek t}} }\else{
{\bigreek t} }\fi\relax}
\def\mattin{\ifmmode{ {\hbox{\smgreek t}} }\else{
{\smgreek t} }\fi\relax}

\begin{document}

\vspace*{0cm}
\begin{center}
{\large \bf
 R\'esultats de compl\'etude pour des classes de types
du syst\`eme ${{\cal AF}2}$}\\
\end{center}

\begin{center}
\bf
Samir FARKH et Karim NOUR\\
\rm
LAMA - Equipe de Logique\\
Universit\'e de Savoie\\
73376 Le Bourget du Lac\\
e-mail sfarkh, knour@univ-savoie.fr\\[1cm]
\end{center}

{\bf Abstract.} J.-L. K\textsc{rivine} introduced the ${{\cal AF}2}$
type system in
order to obtain programs ($\l$-terms) which calculate functions, by
writing
demonstrations of their totalities. We present in this paper two results
of
completness for some types of ${{\cal AF}2}$ and for many notions of
reductions. These results generalize a theorem of R.
L\textsc{abib}-S\textsc{ami} established
in the system ${\cal F}$ of J.-Y. G\textsc{irard}.\\

{\bf R\'esum\'e.} J.-L. K\textsc{rivine} a introduit le syst\`eme de
typage ${{\cal
AF}2}$ pour obtenir des programmes ($\l$-termes) calculant des fonctions
en
\'ecrivant des d\'emonstrations de leur totalit\'e. Nous pr\'esentons
dans
ce papier des r\'esultats de compl\'etude pour certains types de ${{\cal
AF}2}$ et pour plusieurs notions de r\'eductions. Ces r\'esultats
g\'en\'eralisent un th\'eor\`eme de R. L\textsc{abib}-S\textsc{ami}
\'etabli dans le syst\`eme
${\cal F}$ de J.-Y. G\textsc{irard}.\\[1,5cm]

{\Large\bf Introduction} \\

Le syst\`eme de typage ${\cal F}$ a \'et\'e introduit par J.-Y.
G\textsc{irard} (voir
[2]).  Ce syst\`eme est bas\'e sur le calcul propositionnel
intuitionniste
du second ordre, et donc donne la possibilit\'e de quantifier sur les
types. En plus du th\'eor\`eme de normalisation forte qui assure la
terminaison des programmes, le syst\`eme ${\cal F}$ a deux autres
propri\'et\'es :\begin{itemize}
\item[] -- Il permet d'\'ecrire des programmes pour toutes les fonctions
dont la terminaison est d\'emontrable dans l'arithm\'etique de Peano du
second ordre.
\item[] -- Il permet de d\'efinir tous les types de donn\'ees courants :
bool\'eens, entiers, listes, etc.
\end{itemize}
La s\'emantique du syst\`eme ${\cal F}$ propos\'ee par J.-Y.
G\textsc{irard} consiste
\`a associer \`a chaque type $A$ un ensemble de $\l$-termes $\mid A
\mid$,
dans le but d'obtenir le r\'esultat suivant : Si un $\l$-terme $t$ est
de
type $A$, alors il appartient \`a l'ensemble $\mid A \mid$. Ce
r\'esultat
est connu sous le nom du lemme d'ad\'equation, et permet de d\'emontrer
la
normalisation forte du syst\`eme ${\cal F}$.\\
La r\'eciproque du lemme d'ad\'equation n'est pas en g\'en\'eral vraie,
mais R. L\textsc{abib}-S\textsc{ami} a d\'emontr\'e dans [7] un
r\'esultat de ce genre pour
les types clos \`a quantificateurs positifs.\\

 Le syst\`eme de typage ${{\cal AF}2}$, introduit par J.-L.
K\textsc{rivine} (voir [4]), est bas\'e  sur la logique intuitionniste
du second ordre. Au
plan des propri\'et\'es th\'eoriques, normalisation forte et
repr\'esentation des fonctions calculables, le syst\`eme ${{\cal AF}2}$
ne
se distingue pas du syst\`eme ${\cal F}$. La diff\'erence provient de sa
capacit\'e \`a exprimer les sp\'ecifications exactes des programmes, ce
qui
permet d'obtenir un programme calculant une fonction en \'ecrivant une
d\'emonstration du fait que la fonction est du bon type. J.-L.
K\textsc{rivine} a
propos\'e
une s\'emantique pour son syst\`eme, et il a d\'emontr\'e un lemme
d'ad\'equation permettant d'obtenir l'une des plus
importantes propri\'et\'es du  syst\`eme ${{\cal AF}2}$ : C'est
l'unicit\'e
de la repr\'esentation des donn\'ees.\\

Dans ce papier, nous d\'emontrons deux r\'esultats de compl\'etude du
syst\`eme ${{\cal AF}2}$, c'est \`a dire, des \'e\-qui\-va\-lences entre
la
syntaxe du syst\`eme et ses s\'emantiques :
\begin{itemize}
\item Le premier r\'esultat de compl\'etude est obtenu pour les types
\`a
quantificateurs positifs en utilisant une s\'emantique bas\'ee sur les
ensembles des $\l$-termes stables par la $\b \eta$-\'equivalence, ce qui
constitue une g\'en\'eralisation du r\'esultat de R.
L\textsc{abib}-S\textsc{ami} \'etabli
dans le syst\`eme ${\cal F}$.
\item Le second r\'esultat de compl\'etude est \'etabli pour une classe
restreinte de types (les
bons types positifs) qui englobe les types de donn\'ees de J.-L.
K\textsc{rivine} en
utilisant une s\'emantique bas\'ee sur les ensembles des $\l$-termes
stables par la $\b$-expansion. Ce r\'esultat nous permet de comprendre
d'o\`u provient la $\eta$-\'equivalence dans le premier r\'esultat.
\end{itemize}

 Les d\'emonstrations de ces r\'esultats reposent essentiellement sur
des propri\'et\'es syn\-ta\-xiques du syst\`eme ${{\cal AF}2}$ (voir
[9]).\\
 Ces r\'esultats
donnent une r\'eponse partielle \`a la question de la comparaison entre
op\'erateurs de mise en 
m\'emoire s\'emantiques (\`a la K\textsc{rivine}) et syntaxiques (\`a la
N\textsc{our}). Les
deux notions sont en effet identiques dans le cas des classes pour
lesquelles la s\'emantique
de K\textsc{rivine} est compl\`ete (voir [10]).\\

 L'article est organis\'e de la mani\`ere suivante : 

\begin{itemize}
\item[] -- Dans la partie 1, nous rappelons des pr\'eliminaires sur le
$\l$-calcul pur et nous pr\'esentons le syst\`eme de typage
${{\cal AF}2}$ ainsi que ses
propri\'et\'es. Nous donnons, \`a la fin de cette partie, quelques
r\'esultats syntaxiques du syst\`eme que nous utilisons dans les
d\'emonstrations.
\item[] -- La partie 2 est consacr\'ee \`a la s\'emantique
propos\'ee par J.-L. K\textsc{rivine}. Nous rappelons ensuite le lemme
d'ad\'equation
du syst\`eme ${{\cal AF}2}$.
\item[] -- Dans la partie 3, nous d\'emontrons une
r\'eciproque du lemme d'ad\'equation pour les types \`a quantificateurs
du
second ordre positifs avec la $\b \eta$-\'equivalence.
\item[] -- Dans la partie 4, nous pr\'esentons un r\'esultat
analogue pour une classe restreinte de types \`a quantificateurs du
second
ordre positifs (les
bons types), et avec la $\b$-r\'eduction. Nous montrons enfin que les
conditions que nous imposons sur les types sont toutes n\'ecessaires.
\end{itemize}

\section {Le $\l$-calcul pur et typ\'e}

 Nous allons adopter dans cet article les notations de
J.-L. K\textsc{rivine}, par exemple :\\
 On note $\L$ l'ensemble des termes du $\l$-calcul, dits aussi
$\l$-termes.
Etant donn\'es des $\l$-termes $t, u, u_1,..., u_n$,
l'application de $t$ \`a $u$
sera not\'ee $(t)u$, et $(... ((t)u_1) ...)u_n$ sera not\'e
$(t)u_1...u_n$.
On note par
$\f_{\b}$ (resp. $\f_{\eta}$) la $\b$-r\'eduction (resp. la
$\eta$-r\'eduction) et par $\simeq\sb{\b}$ (resp.  $\simeq\sb{\b
\eta}$) la
$\b$-\'equivalence (resp. la $\b \eta$-\'equivalence). Si $t$ est un
$\l$-terme, on note par $Fv (t)$ l'ensemble de ses variables libres.
Alors
on a clairement : $Fv ((t)u) = Fv (t) \cup Fv (u)$ et $Fv (\l xu) = Fv
(u)-
\{x \}$. De plus si $t \f_{\b} t'$, alors $Fv (t') \subseteq Fv (t)$, et
si
$t \f_{\eta}t'$, alors $Fv (t) = Fv (t')$.\\

{\bf Lemme 1.1} {\it Soient $t$ et $t'$ deux $\l$-termes.\\
 Si $t \f_{\b \eta} t'$, alors il existe un $\l$-terme $u$, tel que $t
\f_{\b} u$ et $u \f_{\eta}t'$.}\\

 {\bf Preuve} : Voir [1]. \hfill $\spadesuit$\\

{\bf Lemme 1.2} {\it (i) Si $t$ est $\b \eta$-\'equivalent \`a un terme
normalisable, alors $t$ est normalisable.\\
(ii) Si $t$ est $\b \eta$-\'equivalent \`a un terme clos, alors $t$ est
$\b$-\'equivalent \`a un terme clos.}\\

{\bf Preuve :} Voir [1] et [4]. \hfill $\spadesuit$\\

On consid\`ere le calcul des pr\'edicats intuitionniste du second ordre,
\'ecrit
avec les symboles suivants :
\begin {itemize}
\item Les seuls symboles logiques $\f$ et $\q$ ;
\item Des variables d'individu : $x, y,...$ (appel\'ees aussi variables
du
premier ordre) ;
\item Des variables de relation n-aire ($n=0, 1,...$) : $X,Y,...$
(appel\'ees aussi variables du second ordre) ;
\item Des symboles de fonction n-aire ($n=0, 1,...$) sur les individus ;
\item Des symboles de relation n-aire ($n=0, 1,...$) sur les individus.
\end {itemize}

Chaque variable de relation, et chaque symbole de fonction ou de
relation a
une arit\'e $n \geq 0$ fix\'ee. Un symbole de fonction 0-aire sera
appel\'e
{\bf symbole de constante}. Une variable de relation 0-aire est aussi
appel\'ee {\bf variable propositionnelle}.\\

On suppose qu'il y a  une infinit\'e de variables d'individu, et, pour
chaque $n \geq 0$, une infinit\'e de variables de relation n-aire.\\

La donn\'ee des symboles de fonction et de relation constitue ce qu'on
appelle un {\bf langage}, les autres symboles \'etant communs \`a tous
les
langages.\\

Les {\bf termes} sont construits de la fa\c con suivante :
\begin {itemize}
\item Chaque variable d'individu, et chaque symbole de constante est un
terme ;
\item Si $f$ est un symbole de fonction $n$-aire, et $t_1,...,t_n$ sont
des termes, alors $f(t_1,...,t_n$) est un terme.
\end {itemize}

Les {\bf formules} sont construites de la fa\c con suivante :
\begin {itemize}
\item Si $A$ est une variable ou un symbole de relation $n$-aire, et
$t_1,
..., t_n$ sont des termes, alors $A(t_1,...,t_n$) est une formule dite
{\bf formule atomique} ;

\item Si $A$, $B$ sont des formules, alors $A\f B$ est une formule
;\item
Si $A$ est une formule, alors $\q xA$ et $\q XA$ sont des formules,
$x$ (resp. $X$) \'etant
une variable d'individu (resp. de relation).
\end {itemize}
On d\'efinit les notions de variables {\bf libres} et {\bf li\'ees} de
mani\`ere usuelle.\\

Un terme est dit {\bf clos} s'il n'a pas de variable. Une formule est
dite
{\bf close} si elle
n'a pas de variable libre. La {\bf cl\^oture} d'une formule $F$ est la
formule obtenue en
quantifiant universellement toutes les variables libres de $F$.\\

Soit $\mbox{\boldmath$\xi$} = \xi_1,..., \xi_n$ une suite finie de
variables du premier et/ou du second ordre.
\begin{itemize}
\item La formule $\q \xi_1... \q \xi_nF$ est not\'ee $\q
\mbox{\boldmath$\xi$} F$ ;
\item L'\'ecriture $\ll$ \mbox{\boldmath$\xi$} n'est pas libre dans $A$
$\gg$
signifie que $\xi_i$ ($1 \leq i \leq n$) n'est pas libre dans $A$.
\end{itemize}

La notation $t [u_1/x_1,..., u_n/x_n]$ (resp. $F [u_1/x_1,...,
u_n/x_n]$)
repr\'esente le r\'esultat de la substitution simultan\'ee de $u_1$ \`a
$x_1$,..., $u_n$ \`a $x_n$ dans le terme $t$ (resp. dans la formule
$F$).\\

 Si $X$ est une variable de relation unaire, $t$ et $t'$ deux termes,
alors
la formule $\q X [Xt\f Xt']$ est not\'ee $t=t'$ et est dite {\bf
\'equation
fonctionnelle} ou {\bf formule \'equationnelle}. Un {\bf cas
particulier}
de l'\'equation $t=t'$ est une formule de la forme : \\
 $t [u_1/x_1,..., u_n/x_n] = t' [u_1/x_1,..., u_n/x_n]$ ou $t'
[u_1/x_1,..., u_n/x_n] = t [u_1/x_1,..., u_n/x_n]$,\\
 $u_1,..., u_n$ \'etant des termes du langage.\\

Consid\'erons un langage $L$ du second ordre, et un syst\`eme $E$
 d'\'equations fonctionnelles de $L$. On d\'ecrit un syst\`eme de
$\l$-calcul typ\'e, appel\'e {\bf Arithm\'etique Fonctionnelle du second
 ordre} (en abr\'eg\'e ${{\cal AF}2}$), dont les types sont
les formules de $L$. Dans l'\'ecriture des termes typ\'es
 de ce syst\`eme, nous emploierons les m\^emes symboles pour les
variables
du $\l$-calcul et les variables d'individu du langage $L$. Un {\bf
contexte} $\G$ est un ensemble $x_1 : A_1,..., x_n : A_n$ de
d\'eclarations, o\`u $ x_1,..., x_n$ sont des variables distinctes du
$\l$-calcul, et $A_1,..., A_n$ des formules de $L$.\\

Etant donn\'es un $\l$-terme $t$, un type $A$ et un contexte $\G
= x_1 : A_1,..., x_n : A_n$, on d\'efinit au moyen des r\`egles
suivantes
la notion $\ll$ $t$ est typable, \`a l'aide du syst\`eme \'equationnel
$E$, de
type $A$ dans le contexte $\G$ $\gg$. Cette notion est not\'ee $\G
\v_{{{\cal
AF}2}} t : A$.
\begin{center}
(1) $\G\v_{{{\cal AF}2}} x_i : A_i$ $(1\leq i\leq n)$
\end{center}

\begin{minipage}[t]{250pt}
$  ( 2 ) \quad \displaystyle\frac{ \G, x : A \v_{{\cal AF}2} t : B } {
\G\v_{{\cal AF}2} \l xt : A \f B }$ \\
\end{minipage}
\begin{minipage}[t]{250pt}\sl
$  ( 3 ) \quad  \displaystyle\frac{ \G\v_{{\cal AF}2} u : A \f B \quad
\G\v_{{\cal AF}2}
v : A} { \G\v_{{\cal AF}2} (u)v : B
}$ \\ \end{minipage}

\begin{minipage}[t]{250pt}
$  ( 4 ) \quad \displaystyle\frac{ \G\v_{{\cal AF}2} t : A } {
\G\v_{{\cal
AF}2} t : \q xA }$
{\rm (*)}\\
\end{minipage}
\begin{minipage}[t]{250pt}\sl
$  ( 5 ) \quad  \displaystyle\frac{ \G\v_{{\cal AF}2} t : \q xA} {
\G\v_{{\cal AF}2}
t : A[u/x] }$  {\rm (**)}\\
\end{minipage}

\begin{minipage}[t]{250pt}
$  ( 6 ) \quad \displaystyle\frac{ \G\v_{{\cal AF}2} t : A } {
\G\v_{{\cal
AF}2} t : \q XA }$
{\rm (*)}\\
\end{minipage}
\begin{minipage}[t]{250pt}\sl
$  ( 7 ) \quad  \displaystyle\frac{ \G\v_{{\cal AF}2} t : \q XA} {
\G\v_{{\cal AF}2}
t : A[F/X (x_1,..., x_n)] }$  {\rm (**)}\\
\end{minipage}

\begin{center}
$  ( 8 ) \quad \displaystyle\frac{ \G\v_{{\cal AF}2} t : A[u/x]} {
\G\v_{{\cal AF}2}
t : A[v/x] }$  {\rm (***)}\\
\end{center}

Avec les conditions suivantes : \\

(*) $x$ et $X$ ne sont pas libres dans $\G$.\\
(**) $u$ est un terme et $F$ est une formule. \\
  $A[F/X(x_1,..., x_n)]$ est obtenue en rempla\c cant dans $A$ chaque
formule atomique $X(t_1,...,t_n)$ par $F[t_1/x_1,..., t_n/x_n]$.\\
(***) $u = v$ est un cas particulier d'une \'equation de $E$.\\

{\bf Lemme 1.3} {\it (i) Si $\G\v_{{\cal AF}2} t : A$ et $\G \subseteq
\G'$, alors $\G'\v_{{\cal AF}2} t : A$.\\
(ii) Si $\G\v_{{\cal AF}2} t : A$, alors  $\G'\v_{{\cal AF}2} t : A$,
o\`u
$\G'$ est la restriction de $\G$ aux d\'eclarations contenant les
variables
libres de $t$.}\\

{\bf Preuve} : Par induction sur $\G\v_{{\cal AF}2} t : A$. \hfill
\hfill
$\spadesuit$\\

Le syst\`eme ${{\cal AF}2}$ poss\`ede les propri\'et\'es suivantes :\\

{\bf Th\'eor\`eme 1.4} {\it (i) Si $\G \v_{{\cal AF}2} t : A$, et $t\f
_\beta
t'$, alors $\G \v_{{\cal AF}2} t': A$}.\\
{\it (ii) Si $\G \v_{{\cal AF}2} t : A$, alors $t$ est fortement
normalisable}.\\

{\bf Preuve} : Voir [4]. \hfill \hfill $\spadesuit$\\

Soient $L$ un langage du second ordre et $E$ un syst\`eme d'\'equations
de
$L$. On d\'efinit sur l'ensemble des termes de $L$ une relation
d'\'equivalence not\'ee $\approx_{E}$ de la mani\`ere suivante :
\begin{center}
$a \approx_{E} b \Leftrightarrow  E \v a = b$.
\end{center}
Le lemme suivant d\'ecrit l'\'equivalence qu'on vient de d\'efinir.\\

{\bf Lemme 1.5} {\it $a \approx_{E} b$ ssi on peut l'obtenir au moyen
des
r\`egles suivantes :
\begin{itemize}
\item[] (i) si $a = b$ est un cas particulier d'une \'equation de $E$,
alors $a \approx_{E} b$ ;
\item[] (ii) quels que soient les termes $a, b, c$ de $L$, on a : \\
 -- $a \approx_{E} a$ ; \\
 -- si $a \approx_{E} b$ et $b \approx_{E} c$, alors $a \approx_{E} c$ ;
\item[] (iii) si $f$ est un symbole de fonction $n$-aire de $L$, et si
$a_i
\approx_{E} b_i$ ($1
\leq i \leq n$), alors $f (a_1,..., a_n) \approx_{E} f (b_1,..., b_n)$.
\end{itemize}}

{\bf Preuve} : Voir [4]. \hfill $\spadesuit$ \\

Le lemme suivant permet de g\'en\'eraliser la r\`egle (8) de typage
\'equationnelle.\\

{\bf Lemme 1.6} {\it Si $\G \v_{{\cal AF}2} u : B [a/x]$ et $a
\approx_{E}
b$, alors $\G \v_{{\cal AF}2} u : B [b/x]$.}\\

{\bf Preuve} : Par induction sur la preuve de $a \approx_{E} b$, on
consid\`ere la derni\`ere r\`egle utilis\'ee.
\begin{itemize}
\item Si c'est la r\`egle (i), alors c'est \'evident.
\item Si c'est la r\`egle (ii), alors ou bien $b = a$, et dans ce cas on
a
le r\'esultat, ou bien $a \approx_{E} c$ et $c \approx_{E} b$, et donc
par
hypoth\`ese d'induction, on a $\G \v_{{\cal AF}2} u : B [c/x]$ et $\G
\v_{{\cal AF}2} u : B [b/x]$.
\item Si c'est la r\`egle (iii), alors $a = f (a_1,..., a_n)$ et $b = f
(b_1,..., b_n)$, avec $a_i \approx_{E} b_i$ ($1 \leq i \leq n$). On a
$\G
\v_{{\cal AF}2} u : B [f (a_1,..., a_n)/x]$, donc $\G \v_{{\cal AF}2} u
:
(B [f (x_1,..., a_n)/x]) [a_1/x_1]$, par suite, d'apr\`es l'hypoth\`ese
d'induction, $\G \v_{{\cal AF}2} u : (B [f (x_1,..., a_n)/x])
[b_1/x_1]$.
Ainsi de suite, en utilisant l'hypoth\`ese d'induction $n$ fois, on
obtient
le r\'esultat.  \hfill $\spadesuit$
\end{itemize}

Nous allons pr\'esenter le th\'eor\`eme de programmation pour les
entiers.\\

Etant donn\'es deux $\l$-termes $t, u$
et un entier $k$, on pose, par d\'efinition, $(t)^ku = (t)...(t)u$ (le
$\l$-terme $t$ \'etant
r\'ep\'et\'e $k$ fois au second membre) ; en particulier $(t)^0u = u$.
On
d\'efinit le $\l$-terme $\so{k} = \l x\l f (f)^kx$ ; $\so{k}$ est
appel\'e $\ll$ l'entier $k$ du
$\l$-calcul $\gg$ (ou {\bf entier de Church}).\\

On d\'efinit le {\bf type des entiers naturels} par la formule suivante
:
\begin{center}
$N [x] = \q X\{X0 \f [\q y (Xy\f Xsy)\f  Xx]\}$,
\end{center}
$0$ \'etant une constante pour le z\'ero, et $s$ un symbole de fonction
unaire pour le successeur.\\

Un syst\`eme d'\'equations $E$ est dit {\bf ad\'equat} pour le type des
entiers ssi :
\begin{itemize}
\item[] -- $s (a) \not \approx_{E} 0$ ;
\item[] -- Si $s (a) \approx_{E} s (b)$, alors $a \approx_{E} b$.
\end{itemize}

Soit $f$ une fonction de ${\bf N}^k$ dans ${\bf N}$. On dit que $E$ est
un
syst\`eme d'\'equations {\bf d\'efinissant la fonction} $f$ ssi
:\begin{center}
 $f (s^{n_1} (0),..., s^{n_k} (0)) \approx_{E} s^ {f (n_1,... , n_k)}
(0)$.
\end{center}
Soit $f$ une fonction totale de ${\bf N}^k$ dans ${\bf N}$. Etant
donn\'e
un $\l$-terme $P_f$, on dira que $P_f$ {\bf repr\'esente la fonction}
$f$ si, quels que soient $n_1,..., n_k \in {\bf N}$ :
\begin{center}
$(P_f)\so{n_1}...\so{n_k} \f_{\beta} \so{f (n_1,..., n_k)}$.
\end{center}

{\bf Th\'eor\`eme 1.7 (th\'eor\`eme de programmation)} {\it Soient $f$
une
fonction de ${\bf N}^m$ dans ${\bf N}$, et $E$ un
syst\`eme  d'\'equations ad\'equat d\'efinissant
$f$. Si $P_f$ est un $\l$-terme tel que :
$\v_{{\cal AF}2} P_f : \q x_1...\q x_m \{N [x_1]\f ( ...\f (N [x_m] \f N
[f(x_1,...,
x_m)])...)\}$, \\
alors $P_f$ repr\'esente la fonction $f$}.\\

{\bf Preuve} : Voir [4]. \hfill $\spadesuit$\\

{\bf Exemple.} Consid\'erons par exemple la fonction pr\'ed\'ecesseur 
$p :
{\bf N} \f
{\bf N}$, d\'efinie \`a
l'aide des \'equations : $p (0) = 0$ ; $p (sx) = x$. Un $\l$-terme $t$
qui
repr\'esente
la fonction $p$ (autrement dit un programme pour $p$), est donc un terme
de type (dans le syst\`eme ${{\cal AF}2}$) la formule $\q x \{N [x] \f N
[p
(x)]\}$. \hfill
$\spadesuit$\\

Le th\'eor\`eme de programmation se g\'en\'eralise \`a tous les types de
donn\'ees  syntaxiques de
J.-L. K\textsc{rivine} (voir [4] et [6]).\\

Les r\'esultats que nous allons pr\'esenter maintenant permettent une
simplification dans les d\'emonstrations syntaxiques (voir [9]).\\

On d\'efinit sur les types de ${{\cal AF}2}$ les deux relations binaires
$<'$ et $\sim '$ de la mani\`ere suivante:
\begin {itemize}
\item $\q xA<'A[u/x]$, si $u$ est un terme du langage ;
\item $\q XA<'A[F/X(x_1,..., x_n)]$, si $F$ est une formule du langage ;
\item $A \sim' B$ ssi $A = C[u/x]$, $B = C[v/x]$, et $u = v$ est un cas
particulier d'une \'equation.
\end {itemize}
Soient $\leq$  et $\sim$ les cl\^otures r\'eflexives et transitives
respectives de $<'$ et $\sim '$.\\
On note $A<B$ ssi $A \leq B$ et $A \neq B$.\\

On a le r\'esultat suivant :\\

{\bf Th\'eor\`eme 1.8} {\it $(1)$ Si $\G \v_{{\cal AF}2} x : A$, alors
$A$ s'\'ecrit
$\q \mbox{\boldmath$\xi$} C'$, o\`u \mbox{\boldmath$\xi$} n'est pas
libre dans $\G$, et
 il existe deux types $C$ et $B$ tels que $C\sim C'$, $B \leq C$ et $(x
: B) \in \G$.\\
$(2)$ Si $\G \v_{{\cal AF}2} \l xu : A$, alors $A$ s'\'ecrit
 $\q \mbox{\boldmath$\xi$} (B' \f C')$, o\`u \mbox{\boldmath$\xi$} n'est
pas libre dans $\G$,
 et il existe deux types $C$ et $B$ tels que $C\sim C'$, $B\sim B'$ et 
 $\G, x : B \v_{{\cal AF}2} u : C$.\\
$(3)$ Si $\G \v_{{\cal AF}2} (u) v : A$, alors $A$ s'\'ecrit $\q
\mbox{\boldmath$\xi$} D''$, 
o\`u \mbox{\boldmath$\xi$} n'est pas libre dans $\G$, et
 il existe trois types $F$, $C$ et $D$ tels que $\G \v_{{\cal AF}2} u :
F$, 
 $F \leq C \rightarrow D$,
  $D \leq D'$, $D'\sim D''$ et $\G \v_{{\cal AF}2} v : C$.}\\

{\bf Preuve} : Voir [9]. \hfill $\spadesuit$\\

{\bf Corollaire 1.9} {\it Si $\G,x : A \v_{{\cal AF}2} (x)u_1... u_n :
B$,
alors :\\
$n = 0$, $A\leq C$, $C\sim C'$, $B = \q \mbox{\boldmath$\xi$} C'$, et
\mbox{\boldmath$\xi$}
 n'est libre ni dans $\G$ ni dans $A$,\\
ou\\
$n \not = 0$, $A \leq C_1 \rightarrow  B_1, B'_i \leq C_{i+1}\rightarrow
B_{i+1}$ $(1\leq
i \leq n -1)$, $B'_n\leq B_{n+1}$,
et $B = \q \mbox{\boldmath$\xi$} B'_{n+1}$ o\`u $B_i \sim B'_i$ $(1\leq
i
\leq n +1)$, $\G,x : A \v_{{\cal AF}2} u_i : C_i$
$(1\leq i \leq n)$, et \mbox{\boldmath$\xi$} n'est libre ni dans $\G$ ni
dans $A$}.\\

{\bf Preuve} : Voir [9]. \hfill $\spadesuit$

 \section {$\L$-mod\`ele}
Les s\'emantiques que nous allons pr\'esenter dans ce paragraphe sont
d\^ues \`a
J.-L. K\textsc{rivine} (voir [4]). Elles d\'ependent des classes de
parties de
$\l$-termes appel\'ees
habituellement $\ll$ parties satur\'ees $\gg$. Nous allons utiliser dans
ce papier
deux de ces
s\'emantiques, la premi\`ere est bas\'ee sur les parties stables par la
$\b
\eta$-\'equivalence et
la deuxi\`eme sur les parties stables par la $\b$-expansion.\\

Une partie $G$ de $\L$ est dite {\bf $\f_{\b}$-satur\'ee} (resp. {\bf
$\simeq\sb{\b
\eta}$-satur\'ee}) si, quels que soient les termes $t$ et $u$, on
a : \begin{center}
$u\in G$ et $t \f_{\b} u$ (resp. $t \simeq\sb{\b\eta} u)  \Rightarrow
t \in G$.
\end{center}

 Etant donn\'ees deux parties $G$ et $G'$ de $\L$, on d\'efinit
 une partie de $\L$,
 not\'ee $G\f G'$, en posant : $ u\in (G\f G') \Leftrightarrow (u)t \in
G'$
quel que soit $t\in G$.\\

Il est clair que l'intersection d'un ensemble de parties
$\f_{\b}$-satur\'ees
(resp. $\simeq\sb{\b \eta}$-satur\'ees) de $\L$ est
 $\f_{\b}$-satur\'ee (resp. $\simeq\sb{\b \eta}$-satur\'ee).
De plus si $G'$ est $\f_{\b}$-satur\'ee
(resp. $\simeq\sb{\b \eta}$-satur\'ee), alors $G\f G'$ est
$\f_{\b}$-satur\'ee
 (resp. $\simeq\sb{\b \eta}$-satur\'ee) pour toute
partie $G \subset \L$.\\

Consid\'erons maintenant l'ensemble $P_{S_{\b}(\L)}$ (resp. $P_{S_{\b
\eta}(\L)}$) de toutes les parties $\f_{\b}$-satur\'ees (resp.
$\simeq\sb{\b \eta}$-satur\'ees) de $\L$. Un
sous-ensemble $R$ de $P_{S_{\b}(\L)}$  (resp. $P_{S_{\b \eta}(\L)}$) est
dit {\bf ad\'equat} si :
\begin {itemize}
 \item $G, G' \in R \Rightarrow (G\f G')\in R$ ;
\item Pour toute partie $\s$ de $R$, l'intersection des \'el\'ements de
$\s$ appartient \`a $R$. En particulier $\L \in R$ (prendre $\s =
\emptyset$).
\end {itemize}
 Soit $L$ un langage du second ordre. On va d\'efinir la notion d'un
{\bf
$\L_{\b}$-mod\`ele}
(resp. {\bf $\L_{\b \eta}$-mod\`ele}) pour $L$ :\\

 C'est une modification de la notion classique d'un mod\`ele du second
ordre, dans lequel l'ensemble des valeurs de v\'erit\'e n'est pas
$\{0,1\}$
comme d'habitude, mais un sous-ensemble ad\'equat $R$ de
$P_{S_{\b}(\L)}$
(resp. $P_{S_{\b
\eta}(\L)}$).\\

Un $\L_{\b}$-mod\`ele (resp. $\L_{\b \eta}$-mod\`ele) pour le langage
$L$ est la
donn\'ee de :
 \begin {itemize}
\item Un ensemble $\mid M \mid$ suppos\'e non vide, appel\'e la base de
$M$ ;
 \item Un sous-ensemble ad\'equat $R$ de $P_{S_{\b}(\L)}$ (resp.
$P_{S_{\b
\eta}(\L)}$) ($R$ est appel\'e l'ensemble des  valeurs de v\'erit\'e de
$M$) ;
\item Pour tout symbole de fonction $n$-aire de $L$,
 une fonction $f_M :  \mid M \mid^n \f \mid M \mid$ ;
\item Pour tout symbole de relation $n$-aire $P$ de $L$, une fonction
$P_M
:  \mid M \mid^n \f R$.
\end {itemize}

Soit $M$ un $\L_{\b}$-mod\`ele (resp. $\L_{\b \eta}$-mod\`ele) pour $L$.
{\bf Une interpr\'etation} $I$ est, par d\'e\-fi\-ni\-tion, une
application de
l'ensemble des variables d'individu (resp. de relation $n$-aire) dans
$\mid
M \mid$ (resp. dans $R^{ \mid M \mid^n}$).\\

Soient $I$ une interpr\'etation, $x$ (resp. $X$) une variable d'individu
(resp. de relation $n$-aire), et $a$ (resp. $\Phi$) un \'el\'ement de
$\mid
M \mid$ (resp. de $R^{ \mid M \mid^n}$). On d\'efinit une
interpr\'etation
$J = I [x \leftarrow a]$ (resp. $J = I [X \leftarrow \Phi]$) en posant
$J
(x) = a$ (resp. $J (X) = \Phi$) et $J (\xi) = I (\xi)$ (resp. $J (\xi')
= I
(\xi')$) pour toute variable $\xi \not = x$ (resp. $\xi' \not = X$).\\

 D\'efinissons maintenant la valeur d'une formule de $L$ dans $M$ et
dans une interpr\'etation $I$.\\

 A chaque terme $t$ de $L$ est associ\'ee sa {\bf valeur} $t_{M,I} \in
\mid
M\mid$, d\'efinie par
induction sur $t$ : \begin {itemize}
\item Si $t$ est une variable $x$, alors $t_{M,I} = I (x)$ ;
\item Si $t = f(t^1,..., t^n)$, alors $t_{M,I} = f_M (t_{M,I}^1,...,
t_{M,I}^n)$.
\end {itemize}

Soit $F$ une formule de $L$. La {\bf valeur} de $F$ dans $M$ et dans
l'interpr\'etation $I$, not\'ee par
{\bf $\mid F \mid_{M,I}$} est un \'el\'ement de $R$ d\'efini par
induction
sur $F$ au moyen des
r\`egles suivantes :
\begin{itemize}
\item Si $F$ est une formule atomique $P(t^1,..., t^n)$, o\`u $P$ est un
symbole (resp. une variable) de relation $n$-aire, et $t^1,..., t^n$
sont
des termes de $L$, alors on d\'efinit $\mid F \mid_{M,I}$ par $P_M
(t_{M,I}^1,...,
t_{M,I}^n)$ (resp. $I (X)(t_{M,I}^1,...,t_{M,I}^n))$ qui est un
\'el\'ement
de $R$ ;
\item Si $F$ est $G\f H$, alors $\mid F \mid_{M,I} =  \mid G
\mid_{M,I} \f \mid H \mid_{M,I}$ ;
\item Si $F$ est $\q xG$, o\`u $x$ est une variable d'individu, alors
 $\mid F \mid_{M,I} = \bigcap \{\mid G [x] \mid_{M,I [x \leftarrow a]}$
;
$a\in \mid M \mid\}$ ;
\item Si $F$ est $\q XG$, o\`u $X$ est une variable de relation
$n$-aire,
 alors $\mid F \mid_{M,I} = \\ \bigcap \{ \mid G [X]\mid_{M,I [X
\leftarrow
\Phi]}$ ; $\Phi \in R^{\mid M \mid^n}\}$.
\end {itemize}
 Il est clair que la valeur $\mid F \mid_{M,I}$ ne d\'epend que des
valeurs
dans $I$ des variables libres de $F$. En particulier, si $A$ est un type
clos, $\mid F \mid_{M,I}$ ne d\'epend pas de l'interpr\'etation $I$, on
la
notera alors $\mid F \mid_M$.\\

Soient $M$ un $\L_{\b}$-mod\`ele (resp. $\L_{\b \eta}$-mod\`ele) pour
$L$,
et $u=v$ une \'equation
de $L$. On dit que {\bf $M$ satisfait $u = v$}, si la cl\^oture de cette
formule est vraie dans $M$.\\

Si {\bf $E$} est un ensemble d'\'equations de $L$, on dit que $M$ {\bf
satisfait} $E$, ou que $M$ est un {\bf mod\`ele de} $E$, si $M$
satisfait
toute \'equation de $E$.\\

Pour tout type clos $A$, on note par $\mid A \mid_{\b} = \bigcap \{ \mid
A
\mid_M$ / $M$ $\L_{\b}$-mod\`ele de $E$ \}, et par $\mid A \mid_{\b
\eta} = \bigcap \{ \mid A \mid_M$ / $M$ $\L_{\b \eta}$-mod\`ele de
$E$ \}.\\

{\bf Th\'eor\`eme 2.1 (lemme d'ad\'equation)} {\it Soient $E$ un
ensemble
fini d'\'equations d'un langage $L$, $t$ un $\l$-terme, et $A$ un type
clos
du syst\`eme ${{\cal AF}2}$.\\
Si $\v_{{\cal AF}2} t : A$, alors $t \in \mid A \mid_{\b}$ (resp. $t \in
\mid A \mid_{\b \eta}$) }.\\

{\bf Preuve} : Voir [5]. \hfill $\spadesuit$\\

Dans la d\'emonstration de ce th\'eor\`eme, on utilise les deux lemmes
suivants :\\

{\bf Lemme 2.2} {\it Soient $t$ un terme, et $A$ une formule de $L$
ayant
$x$ comme
variable libre, alors : $\mid A[t/x] \mid_{M,I} = \mid A \mid_{M,I [x
\leftarrow t_{M,I}]}$}.\\

{\bf Lemme 2.3} {\it Soit $F$ une formule de $L$, ayant $x_1,..., x_n$
comme variables
libres, et soit $\Phi \in
R^{ \mid M \mid^n}$ d\'efinie par $\Phi (a_1,..., a_n) = \mid
F \mid_{M,I [x_1\leftarrow a_1,..., x_n\leftarrow a_n]}$
pour tout $a_1,..., a_n \in \mid M \mid$. Si $A$ est une formule ayant
$X$
(variable de relation
$n$-aire) comme variable libre, alors : $\mid A[F/X (x_1,..., x_n)]
\mid_{M,I} = \mid A \mid_{M,I [X \leftarrow \Phi]}$}.\\

{\bf Remarque}. Pour avoir le lemme d'ad\'equation, il suffit de donner
la
d\'efinition suivante d'ensemble satur\'e $G$ : quels que soient les
termes $t,t_1,..., t_n, u$, \\ $(u[t/x])t_1...t_n\in G \Rightarrow (\l
xu)tt_1...t_n\in G$.\\

Dans la suite nous allons d\'emontrer des r\'eciproques du th\'eor\`eme
2.1.

\section {Premier r\'esultat de compl\'etude}

Dans ce paragraphe, nous allons g\'en\'eraliser le r\'esultat de
R. L\textsc{abib}-S\textsc{ami} pour le syst\`eme ${{\cal AF}2}$, et
pour une classe plus
large de types.\\

On d\'efinit de la fa\c con suivante les types {\bf \`a quantificateurs
positifs} (resp. {\bf \`a quantificateurs n\'egatifs}), not\'es en
abr\'eg\'e  $\q_2^+$ (resp. $\q_2^-$) :
 \begin{itemize}
\item Une formule atomique est  $\q_2^+$ (resp. $\q_2^-$) ;
\item Si $A$ est $\q_2^+$ (resp. $\q_2^-$) et $B$ est $\q_2^-$ (resp.
$\q_2^+$), alors $B \rightarrow A$ est $\q_2^+$ (resp. $\q_2^-$) ;
\item Si $A$ est $\q_2^+$, et $x$ (resp. $X$) une variable d'individu
(resp. de relation $n$-aire), alors $\q xA$ (resp. $\q XA$) est $\q_2^+$
;
\item Si $A$ est $\q_2^-$ et $x$ (resp. $X$) une variable d'individu
(resp.
de relation $n$-aire qui ne figure pas dans $A$), alors $\q xA$ (resp.
$\q
XA$) est $\q_2^-$.
\end{itemize}
Les formules $\q_2^+$ (resp. $\q_2^-$) sont les formules o\`u les
quantificateurs du second ordre $\ll$ ac\-tifs $\gg$ sont en position
positive (resp.
n\'egative) dans la formule.\\

Il est facile de voir que : Si $A$ est $\q_2^+$ (resp. $\q_2^-$) et $A
\sim
B$, alors $B$ est
$\q_2^+$ (resp. $\q_2^-$). De plus si $A$ est $\q_2^-$ et $A \leq B \f
C$, alors $B$ est $\q_2^+$
et $C$ est $\q_2^-$.\\

K. N\textsc{our} a d\'efini dans [8] le syst\`eme de typage ${{\cal
AF}2_0}$ qui
n'est autre que le syst\`eme ${{\cal AF}2}$ o\`u on remplace la r\`egle
de
typage (7) par la r\`egle :

\begin{center}
$(7_0) \quad  \displaystyle\frac{ \G\v_{{\cal AF}2_0} t : \q XA}
{\G\v_{{\cal AF}2_0} t : A[Y/X (x_1,..., x_n)] }$\\
\end{center}
o\`u $Y$ est une variable ou symbole de relation de m\^eme arit\'e que
$X$.\\

{\bf Th\'eor\`eme 3.1} {\it Soient A un type $\q_2^+$ du syst\`eme
${{\cal
AF}2}$, et $t$ un
$\l$-terme normal clos.\\
Si $\v_{{\cal AF}2} t : A$, alors $\v_{{\cal AF}2_0} t : A$.}\\

{\bf Preuve} : Voir [8]. \hfill $\spadesuit$\\

Soient $L$ un langage du second ordre et $E$ un syst\`eme d'\'equations
de
$L$, ce qui d\'efinit le
syst\`eme de typage ${{\cal AF}2}$.\\

On se propose de d\'emontrer le th\'eor\`eme suivant :\\

{\bf Th\'eor\`eme I  Soient A un type $\q_2^+$ clos du syst\`eme ${{\cal
AF}2}$, et $t$ un
$\l$-terme, alors :\\ $t \in \mid A \mid_{\b \eta}$ $\so{ssi}$ il existe
un $\l$-terme $t'$
 tel que $t\simeq\sb{\b \eta} t'$ et $\v_{{\cal AF}2} t' : A$.}\\

Pour d\'emontrer ce th\'eor\`eme, on a besoin de certaines d\'efinitions
et
d'un lemme.\\

Soient $\Omega = \{ x_i/i \in {\bf N} \}$ une \'enum\'eration d'un
ensemble
infini
de variables du
$\l$-calcul, et $\{A_i/i \in
{\bf N} \}$ une \'enum\'eration des types $\q_2^-$ de ${{\cal AF}2}$, 
o\`u
chaque type $\q_2^-$
se r\'ep\`ete une infinit\'e de fois. On d\'efinit alors l'ensemble
$\mit
\G^- = \{x_i : A_i /i \in {\bf N} \}$. Soit $u$ un $\l$-terme, tel que
$Fv
(u) \subseteq \Omega$. On d\'efinit le contexte $\mit \G^-_u$ comme
\'etant
la restriction de $\mit \G^-$ sur les d\'eclarations contenant les
variables de $Fv (u)$. La notation $\mit \G^- \v_{{\cal AF}2} u : B$
exprime que $\mit \G^-_u \v_{{\cal AF}2} u : B$.\\

On pose $\mit \G^- \v_{{\cal AF}2}^{\b \eta} u : B$ ssi il existe un
$\l$-terme $u'$, tel que $u \simeq\sb{\b \eta} u'$ et $\mit \G^-
\v_{{\cal
AF}2} u' : B$.\\

Consid\'erons $M_0$ l'ensemble de tous les termes de $L$. On d\'efinit
un
$\L_{\b \eta}$-mod\`ele ${\cal M}_{\G^-}$ (not\'e dans la suite ${\cal
M}$) de la fa\c con suivante :
\begin {itemize}
\item L'ensemble de base est ${\cal |M|} =
\displaystyle\frac{M_0}{\approx_{E}}$ (l'ensemble des classes
d'\'equivalence modulo la relation $\approx_{E}$) ;
\item L'ensemble ad\'equat $R = P_{S_{\b \eta}}(\L)$ ;
\item L'interpr\'etation de chaque symbole de fonction $n$-aire $f$ est
l'application $f_{\cal M}$ de $\mid {\cal M} \mid ^n$ dans $\mid {\cal
M}
\mid$ d\'efinie par $f_{\cal M} (\sou{a_1},..., \sou{a_n}) = \sou{f
(a_1,..., a_n)}$ ;
\item L'interpr\'etation de chaque symbole de relation $n$-aire $P$ est
l'application $P_{\cal M}$
de\\ $\mid {\cal M} \mid ^n$ dans $P_{S_{\b \eta}}(\L)$ d\'efinie par
$P_{\cal M} (\sou{a_1},..., \sou{a_n})
= \{ \t \in \L : \mit \G^- \v_{{\cal AF}2}^{\b \eta} \t : P(a_1,...,
a_n)\}$.
\end {itemize}

Ensuite on d\'efinit une interpr\'etation ${\cal I}$ sur les variables
en
posant :
\begin {itemize}
\item ${\cal I} (x) = \sou{x}$, o\`u $\sou{x}$ est la classe de $x$
modulo
$\approx_{E}$ ;
\item ${\cal I} (X) = \Phi$, o\`u $\Phi$ est
l'application de $\mid {\cal M} \mid^n$ dans $P_{S_{\b \eta}}(\L)$
d\'efinie par
 $\Phi (\sou{a_1},..., \sou{a_n}) = \{ \t \in \L : \mit \G^- \v_{{\cal
AF}2}^{\b \eta} \t : X(a_1,
..., a_n)\}$.
\end {itemize}

 Les $f_{\cal M}, P_{\cal M}$, et $\Phi$ sont bien d\'efinies. En effet
:\\

Si $(\sou{a_1},..., \sou{a_n}) = (\sou{b_1},..., \sou{b_n})$, alors $a_i
\approx_{E} b_i$ ($1 \leq i \leq n$). Donc, d'apr\`es le lemme 1.5, $f
(a_1,..., a_n) \approx_{E} f (b_1,..., b_n)$, et par cons\'equent
$f_{\cal
M} (\sou{a_1},..., \sou{a_n}) = f_{\cal M} (\sou{b_1},...,
\sou{b_n})$.\\

De m\^eme supposons que $(\sou{a_1},..., \sou{a_n}) = (\sou{b_1},...,
\sou{b_n})$, donc $a_i \approx_{E} b_i$ ($1 \leq i \leq n$). D'o\`u,
d'apr\`es le lemme 1.6, $\{ \t \in \L : \mit \G^- \v_{{\cal AF}2}^{\b
\eta}
\t : P(a_1,..., a_n)\} = \{ \t \in \L : \mit \G^- \v_{{\cal AF}2}^{\b
\eta}
\t : P(b_1,..., b_n)\}$, ou alors, $P_{\cal M} (\sou{a_1},...,
\sou{a_n}) =
P_{\cal M} (\sou{b_1},..., \sou{b_n})$.\\

La m\^eme d\'emonstration se fait pour $\Phi$. \hfill $\spadesuit$\\

Le lemme suivant va nous permettre de d\'emontrer le th\'eor\`eme I.\\

{\bf Lemme 3.2} {\it Soient $S$ une formule du langage $L$, et $\t$ un
 $\l$-terme.\\
(i) Si $S$ est $\q_2^+$, et $\t \in \mid S\mid_{\cal {M,I}}$, alors
$\mit
\G^-\v_{{\cal AF}2}^{\b \eta} \t : S$.\\
(ii) Si $S$ est $\q_2^-$, et $\mit \G^-\v_{{\cal AF}2}^{\b \eta} \t :
S$,
alors $\t \in \mid S\mid_{\cal {M,I}}$}.\\

{\bf Preuve} : Par induction simultan\'ee sur les types $\q_2^+$ et
$\q_2^-$.\\

\so{Preuve de (i)}
\begin{itemize}
\item Si $S$ est atomique, alors $S = X(a_1,...,a_n)$,  o\`u les $a_i$
sont
des termes, et $X$ une variable (resp. un symbole) de relation $n$-aire.
Soit $\t \in \mid S \mid_{\cal {M,I}} = \mid X (a_1,..., a_n) \mid_{\cal
{M,I}} = \{\th \in \L : \mit \G^- \v_{{\cal AF}2}^{\b \eta} \th : X
(a_1,..., a_n)\}$. Il en r\'esulte que $\t \in \mid S \mid_{\cal {M,I}}$
ssi $\mit \G^- \v_{{\cal AF}2}^{\b \eta} \t : S$.
\item Si $S = B\f C$, o\`u $B$ est $\q_2^-$ et $C$ est $\q_2^+$,  soit
$\t \in \mid B\f C \mid_{\cal {M,I}}$. Il existe une infinit\'e de $i
\in {\bf
N}$ tel que $B = A_i$. On choisit un $i$ de fa\c con que $x_i$ ne soit
pas
libre dans $\t$. On a $x_i : B \v_{{\cal AF}2} x_i : B$, donc, d'apr\`es
(ii), $x_i \in \mid B \mid_{\cal {M,I}}$, et alors $(\t)x_i\in \mid C
\mid_{\cal {M,I}}$, et donc, d'apr\`es l'hypoth\`ese d'induction, $\mit
\G^- \v_{{\cal AF}2}^{\b \eta}
(\t)x_i : C$. Donc il existe un $\l$-terme $\t'$ tel que $(\t)x_i
\simeq\sb{\b \eta} \t'$ et $\mit \G^-_{\t'}
\v_{{\cal AF}2} \t': C$, par suite $\t \simeq\sb{\b \eta} \l x_i(\t)x_i
\simeq\sb{\b \eta} \l x_i \t'$. Si $x_i \in Fv (\t')$, alors $\mit
\G^-_{\l
x_i \t'} \v_{{\cal AF}2} \l x_i \t' : B\f C$. Sinon on peut \'ecrire
$\mit
\G^-_{\t'}, x_i : B \v_{{\cal AF}2} \t': C$, et comme $Fv (\t') = Fv (\l
x_i \t')$, alors $\mit \G^-_{\l x_i \t'} \v_{{\cal AF}2} \l x_i \t' :
B\f
C$. Par cons\'equent $\mit \G^- \v_{{\cal AF}2}^{\b \eta} \t : S$.
\item Si $S=\q XB$, avec $B$ est $\q_2^+$, alors $\t\in \mid \q XB
\mid_{\cal {M,I}} \Leftrightarrow (\q \Phi : \mid {\cal M} \mid^n \f
P_{S_{\b \eta}(\L)}) (\t\in \mid B[X] \mid_{{{\cal M,I} [X \leftarrow
\Phi]}}$.
 Soit $Y$ une variable de relation $n$-aire qui ne figure pas dans $\mit
\G^-_{\t}$. Donc $\t\in \mid B[X] \mid_{{{\cal M,I} [X \leftarrow \mid Y
\mid_{\cal {M,I}}]}} = \mid B[Y / X] \mid_{\cal {M,I}}$, d'apr\`es le
lemme
2.3. D'o\`u par hypoth\`ese d'induction, $\mit \G^- \v_{{\cal AF}2}^{\b
\eta} \t : B[Y]$, donc il existe un $\l$-terme $\t'$ tel que $\t
\simeq\sb{\b \eta} \t'$ et $\mit \G^-_{\t'}
\v_{{\cal AF}2} \t' : B[Y]$. Par cons\'equent il existe un $\l$-terme
$u$
tel que $\t \f_{\b \eta} u$ et $\t' \f_{\b \eta} u$, et, d'apr\`es le
lemme
1.1, il existe un $\l$-terme $\t''$ tel que $\t'\f_{\b}\t''$ et
$\t''\f_{\eta}u$. Donc $Fv (u) = Fv (\t'') \subseteq Fv (\t')$ et $\mit
\G^-_{\t''} \v_{{\cal AF}2} \t'' : B[Y]$ , d'o\`u, par le choix de $Y$
et
le fait que $Fv (\t'') \subseteq Fv (\t)$, on d\'eduit que $\mit
\G^-_{\t''} \v_{{\cal AF}2} \t'' : \q YB[Y]$, et donc $\mit \G^-
\v_{{\cal
AF}2}^{\b \eta} \t : S$.
\item Si $S = \q xB$, avec $B$ est $\q_2^+$ , alors $\t \in \mid \q xB
\mid_{\cal {M,I}} \Leftrightarrow (\q a \in \mid {\cal M} \mid ) (\t \in
\mid
B[x] \mid_{{{\cal M,I} [x \leftarrow a]}})$. Soit $y$ une variable de
$L$
qui ne figure pas dans les formules de $\mit \G^-_{\t}$. Alors on a :
$\sou{y} \in \mid {\cal M} \mid$, donc $\t \in \mid B[x] \mid_{{{\cal
M,I}
[x \leftarrow \sou{y}]}} = \mid
B[y/x] \mid_{\cal {M,I}}$ (d'apr\`es le lemme 2.2). D'o\`u, d'apr\`es
l'hypoth\`ese d'induction, $\mit \G^- \v_{{\cal AF}2}^{\b \eta} \t :
B[y]$, et
par le m\^eme raisonnement qu' au cas pr\'ec\'edent et le choix de $y$,
on obtient
$\mit \G^- \v_{{\cal AF}2}^{\b \eta} \t : S$.
\end {itemize}

\so{Preuve de (ii)}
{\begin {itemize}
\item Si $S$ est atomique, le r\'esultat d\'ecoule imm\'ediatement de la
d\'efinition de $I_{\cal M}$ (voir (i)).
\item Si $S = B \rightarrow C$, o\`u $B$ est $\q_2^+$ et $C$ est
$\q_2^-$,
supposons $\mit \G^- \v_{{\cal AF}2}^{\b \eta} \t : B \rightarrow C$.
Donc
 il existe un $\l$-terme $\t'$ tel que
$\t \simeq\sb{\b \eta} \t'$ et $\mit \G^-_{\t'}\v_{{\cal AF}2} \t' : B
\rightarrow C$. Si $u \in \mid B \mid_{\cal {M,I}}$,
alors d'apr\`es (i), $\mit \G^-\v_{{\cal AF}2}^{\b \eta}  u : B$, donc
 il existe un $\l$-terme $u'$ tel que $u
\simeq\sb{\b \eta} u'$ et $\mit \G^-_{u'}\v_{{\cal AF}2} u' : B$. D'o\`u
$\mit \G^-_{(\t')u'}\v_{{\cal AF}2} (\t') u' : C$, et comme $(\t)u
\simeq\sb{\b \eta} (\t')u'$, alors $\mit \G^-\v_{{\cal AF}2}^{\b \eta}
(\t)u : C$. D'o\`u, d'apr\`es l'hypoth\`ese d'induction, $(\t)u \in \mid
C
\mid_{\cal {M,I}}$. Par cons\'equent $\t \in \mid B \f C \mid_{\cal
{M,I}}$.
\item Si $S = \q xB$, o\`u $B$ est  $\q_2^-$, supposons que $\mit
\G^-\v_{{\cal AF}2}^{\b \eta} \t : \q xB$ (-), et soit $a \in \mid {\cal
M}
\mid$. Donc $a = \sou{b}$, o\`u $b$ est un terme de $L$. D'apr\`es (-),
il existe
 un $\l$-terme $\t'$ tel que $\t
\simeq\sb{\b \eta} \t'$ et $\mit \G^-_{\t'}\v_{{\cal AF}2} \t' : \q xB$,
par suite $\mit \G^-_{\t'}\v_{{\cal AF}2} \t' : B [b/x]$. Donc,
d'apr\`es
l'hypoth\`ese
d'induction, $\t' \in \mid B[b/x] \mid_{\cal {M,I}} = \mid B
\mid_{{{\cal
M,I} [x \leftarrow \sou{b}]}} = \mid B \mid_{{{\cal M,I} [x \leftarrow
a]}}$ (d'apr\`es le lemme 2.2). Comme $\mid B \mid_{{{\cal M,I} [x
\leftarrow
a]}}$ est $\simeq\sb{\b \eta}$-satur\'ee, on aura $\t \in \mid B
\mid_{{{\cal M,I} [x \leftarrow a]}}$.
\item Si $S = \q XB$, o\`u $B$ est  $\q_2^-$, supposons que $\mit
\G^-\v_{{\cal AF}2}^{\b \eta} \t : \q XB$. Alors il existe un $\l$-terme
$\t'$ tel que 
$\t \simeq\sb{\b \eta}
\t'$ et $\mit \G^-_{\t'}\v_{{\cal AF}2} \t' : B$, donc $\mit \G^-
\v_{{\cal
AF}2} \t' : B$, et par hypoth\`ese d'induction, $\t' \in \mid B
\mid_{\cal
{M,I}}$. Donc $\t' \in \mid \q X B \mid_{\cal {M,I}}$, car $X$ ne figure
pas dans $B$, et comme $\mid \q X B \mid_{\cal {M,I}}$ est $\simeq\sb{\b
\eta}$-satur\'ee,
 on aura $\t \in \mid \q XB \mid_{\cal {M,I}}$.  \hfill $\spadesuit$
\end {itemize}

On peut d\'eduire maintenant la preuve du th\'eor\`eme I.\\

{\bf Preuve du th\'eor\`eme I}
\begin{itemize}
\item La condition suffisante est une cons\'equence directe du
th\'eor\`eme 2.1,
et du fait que l'interpr\'etation d'un type du syst\`eme ${{\cal AF}2}$
dans un
$\L_{\b \eta}$-mod\`ele $M$ est une partie $\simeq\sb{\b
\eta}$-satur\'ee.
\item La condition est n\'ecessaire. En effet : Soient $A$ un type
$\q_2^+$
clos, et
$t$ un $\l$-terme tel que $t \in \mid A \mid_{\b \eta}$, alors $t \in
\mid A \mid_{\cal M}$, pour tout $\L_{\b \eta}$-mod\`ele ${\cal M}$
associ\'e \`a un ensemble $\mit \G^-$ (comme d\'ecrit avant). De plus on
peut supposer que $\mit \G^-$ ne contient pas de d\'eclarations pour les
variables libres de $t$, i.e $\mit \G^-_t = \emptyset$. D'apr\`es le (i)
du
lemme 3.2, $\mit \G^-\v_{{\cal AF}2}^{\b \eta} t : A$, donc il existe un
$\l$-terme
$t'$ tel que $t \simeq\sb{\b
\eta} t'$ et $\mit \G^-_{t'}\v_{{\cal AF}2} t' : A$, par suite, il
existe
un $\l$-terme $u$ tel que $t \f_{\b \eta} u$ et $t' \f_{\b \eta} u$, et
d'apr\`es le lemme 1.1, il existe un $\l$-terme $t''$ tel que $t'
\f_{\b}
t''$ et $t'' \f_{\eta} u$. D'o\`u $Fv (t'') \subseteq Fv (t')$ et $\mit
\G^-_{t''} \v_{{\cal AF}2} t'' : A$. Or $Fv (t'') = Fv (u) \subseteq Fv
(t)$, donc $\mit \G^-_{t''} \subseteq \mit \G^-_t = \emptyset$, par
cons\'equent $t \simeq\sb{\b \eta} t''$ et $\v_{{\cal AF}2} t'' : A$.
\hfill $\spadesuit$\\
\end {itemize}

 D'apr\`es le th\'eor\`eme 3.1, on peut d\'eduire le r\'esultat suivant
:\\

{\bf Th\'eor\`eme 3.3} {\it Soient A un type $\q_2^+$ clos du syst\`eme
${{\cal AF}2}$, et $t$ un
$\l$-terme, alors :\\ $t \in \mid A \mid_{\b \eta}$ $\so{ssi}$ il existe
un $\l$-terme $t'$
tel que $t\simeq\sb{\b \eta} t'$ et $\v_{{\cal AF}2_0} t' : A$.}\\

{\bf Corollaire 3.4} {\it Soient $A$ un type $\q_2^+$ clos du syst\`eme
${{\cal AF}2}$, et $t$ un
$\l$-terme.\\
Si $t \in \mid A \mid_{\b \eta}$, alors $t$ est normalisable et
$\b$-\'equivalent \`a un terme clos.}\\

{\bf Preuve} : Si $t \in \mid A \mid_{\b \eta}$, alors, d'apr\`es le
th\'eor\`eme I, il existe un $\l$-terme $t'$ tel que $t \simeq\sb{\b
\eta} t'$
et $\v_{{\cal AF}2} t' : A$. Donc $t$ est $\b\eta$-\'equivalent \`a un
terme normalisable clos. D'o\`u le r\'esultat d'apr\`es le lemme 1.2.
\hfill $\spadesuit$\\

{\bf Remarque.} La condition $\q_2^+$ est n\'ecessaire pour avoir le
th\'eor\`eme I. En effet,
soit
\begin{center}
$D = \q X\{\q Y (Y\f X)\f X\}$
\end{center}
Il est clair que $D$ n'est pas $\q_2^+$. Posons  $t = \l x (x)
(\delta)\delta$, o\`u $\delta = \l x (x)x$, alors on a: $t \in \mid D
\mid_{\b \eta}$, et $t$ n'est pas normalisable. Pour montrer que $t$ est
bien un terme de
$\mid D \mid_{\b \eta}$, soit $M$ un
$\L_{\b \eta}$-mod\`ele, montrons que $t \in \mid \q Y (Y \f X)\f X
\mid_{M, [X \leftarrow \Xi]}$  pour tout $\Xi \in P_{S_{\b \eta}}(\L)$.
Soient $u \in \mid \q Y (Y
\f X)\mid_{M, [X \leftarrow \Xi]}$, et $\Xi'$ l'ensemble des $\l$-termes
qui ne sont pas
normalisables. $\Xi'$ est
\'evidemment une partie $\simeq\sb{\b \eta}$-satur\'ee de $\L$.
 Comme $u \in \mid \q Y (Y \f X)\mid_{M, [X \leftarrow \Xi]}$, alors
 $u \in \mid Y \f X \mid_{M, [X \leftarrow \Xi, Y \leftarrow \Xi']}$,
c'est \`a dire
 $u \in \Xi' \f \Xi$. On a $(\delta)\delta \in \Xi'$, donc
$(u)(\delta)\delta \in \Xi$. Or $(t)u = (\l
x (x) (\delta)\delta)u \f_{\b} (u)(\delta)\delta$, d'o\`u $(t)u \in
\Xi$.\\
Posons  $t' = \l x (x) y$, o\`u $y$ est une variable. On a : $t' \in
\mid D
\mid_{\b \eta}$, et $t'$ n'est pas $\b\eta$-\'equivalent \`a un terme
clos.
En effet, d'une part $t'$ est normal et non clos, d'autre part nous
montrons
 que $t' \in \mid D \mid_{\b \eta}$, soient $M$
un $\L_{\b \eta}$-mod\`ele quelconque, et $u \in \mid \q Y (Y \f
X)\mid_{M, [X \leftarrow \Xi]}$
 pour tout $\Xi \in P_{S_{\b \eta}}(\L)$. Consid\'erons $\Xi''$
l'ensemble des $\l$-termes qui sont $\b \eta$-\'equivalent \`a une
variable. $\Xi''$ est \'evidemment une partie $\simeq\sb{\b
\eta}$-satur\'ee de
$\L$. Alors $u \in \mid Y \f X \mid_{M, [X \leftarrow \Xi, Y \leftarrow
\Xi'']}$, c'est \`a dire
 $u \in \Xi'' \f \Xi$, et comme $y \in \Xi''$, on obtient $(u)y \in
\Xi$. Or
$(t')u = (\l x (x) y)u \f_{\b} (u)y$, d'o\`u $(t')u \in \Xi$. \hfill
$\spadesuit$ \\

Le syst\`eme de typage ${\cal F}$ de J.-Y. G\textsc{irard} est le
sous-syst\`eme de
${{\cal AF}2}$, o\`u on a seulement des variables propositionnelles et
des
constantes (symboles de relation $0$-aire). Donc les variables du
premier
ordre, les symboles de fonction et le syst\`eme d'\'equations sont
inutiles. Les r\`egles de typages sont les r\`egles $(1), (2), (3)$ et
$(6), (7)$ du syst\`eme
${{\cal AF}2}$ restreintes aux variables propositionnelles.\\

Si on se restreint au syst\`eme ${\cal F}$, un $\L_{\b}$-mod\`ele (resp.
$\L_{\b \eta}$-mod\`ele)
est compos\'e uniquement d'une partie ad\'equate $R$ de $P_{S_{\b}(\L)}$
(resp. $P_{S_{\b
\eta}(\L)}$), et pour toute constante (resp. variable) propositionnelle
$P$, d'un \'el\'ement
$\mid P \mid_M$ de $R$.\\

 Le th\'eor\`eme I s'\'enonce dans le syst\`eme ${\cal F}$ de la fa\c
con
suivante :\\

{\bf Th\'eor\`eme 3.5} {\it Soient $A$ un type $\q_2^+$ clos du
syst\`eme
${\cal F}$, et $t$ un
$\l$-terme.\\
$t \in \mid A \mid_{\b \eta}$ $\so{ssi}$ il existe un $\l$-terme $t'$
tel que
 $t \simeq\sb{\b\eta} t'$ et $\v_{\cal F} t' : A$.}\\

Ce th\'eor\`eme a \'et\'e d\'emontr\'e par R.
L\textsc{abib}-S\textsc{ami} (voir [7]).

\section {Second r\'esultat de compl\'etude}

Dans ce paragraphe, nous allons d\'emontrer un th\'eor\`eme analogue au
th\'eor\`eme I avec la
${\b}$-r\'eduction et pour une classe plus restreinte de types.\\

Les types {\bf propres} sont d\'efinis de la fa\c con suivante :
\begin {itemize}
\item Une formule atomique est propre ;
\item Si $A$ et $B$ sont propres, alors $A \rightarrow B$ est propre ;
\item Si $A$ est propre, et $x$ (resp. $X$) une variable d'individu
(resp.
de relation $n$-aire qui figure dans $A$), alors $\q xA$ (resp. $\q XA$)
est propre.
\end {itemize}

On d\'efinit de la fa\c con suivante les types {\bf positifs} (resp.
{\bf
n\'egatifs}), not\'es en abr\'eg\'e  $\q^+$ (resp. $\q^-$) :

\begin {itemize}
\item Une formule atomique est  $\q^+$ (resp. $\q^-$) ;
\item Si $A$ est $\q^+$ (resp. $\q^-$ qui ne commence pas par un
quantificateur du premier ordre)
 et $B$ est $\q^-$ (resp. $\q^+$), alors $B \rightarrow A$ est $\q^+$
(resp. $\q^-$) ;
\item Si $A$ est $\q^+$, et $X$ une variable de relation qui figure dans
$A$, alors $\q XA$ est $\q^+$ ;
\item Si $A$ est $\q^-$, alors $\q xA$ est $\q^-$.
\end {itemize}

D'apr\`es cette d\'efinition chaque type $\q^+$ (resp. $\q^-$) est
propre. \\

Pr\'ecis\'ement les types $\q^+$ sont les types
o\`u :
\begin{itemize}
\item[] -- Les quantificateurs du second ordre $\ll$ actifs $\gg$ sont
positifs.
\item[] -- Les quantificateurs du premier ordre sont n\'egatifs.
\item[] -- Les seules variables du second ordre sur lesquelles on peut
quantifier
 sont $\ll$ actives $\gg$.
\item[] -- On n'a pas le droit de mettre un quantificateur du premier
ordre
juste derri\`ere une
fl\`eche.
\end{itemize}
D'apr\`es la d\'efinition ci-dessus, on peut remarquer que si $A$ est
$\q^-$, alors $A$ est de la forme  $\q \mbox{\boldmath$\xi$} (A_1 \f
(A_2\f (...\f (A_n \f X
(t_1,..., t_n))...)))$, o\`u \mbox{\boldmath$\xi$} est une suite finie
de
variables du premier ordre, $A_i$ ($1 \leq i \leq n)$ des types $\q^+$,
et
$X$ une variable ou un symbole de relation $n$-aire.\\

Il est clair que si $A$ est un type $\q^-$ et $A \leq B$, alors $B$ est
$\q^-$, et une variable du second ordre qui est libre dans $B$ est aussi
libre dans $A$.\\

Soit $A$ un type. On d\'efinit les {\bf sous-types positifs} (resp. {\bf
n\'egatifs}) 
de $A$ de la fa\c con inductive suivante :
\begin{itemize}
\item[] -- Si $A = X (t_1,..., t_n)$, alors $A$ est le seul sous-type
positif et n\'egatif de $A$.
\item[] -- Si $A = B \f C$, alors les sous-types positifs (resp.
n\'egatifs) de $A$
sont les sous-types n\'egatifs (resp. positifs) de $B$, et les
sous-types positifs
 (resp. n\'egatifs) de $C$.
\item[] -- Si $A = \q XB$ (resp. $\q xB$), alors les sous-types positifs
 (resp. n\'egatifs) de $A$ sont les sous-types positifs (resp.
n\'egatifs) de $B$. 
 \end{itemize}
 Les sous-types positifs (resp. n\'egatifs) d'un type forment en
g\'en\'eral un ensemble infini, car par exemple pour le
type $\q XB [X]$, il faut consid\'erer comme sous-types positifs, $B
[Y]$, $B
[Z]$, etc. Intuitivement, on \'ecrit la formule avec des noms
diff\'erents
pour les variables li\'ees auquelles on associe des ensembles
d\'enombrables disjoints de variables. Pour chacune des variables
li\'ees,
on utilise l'ensemble correspondant de variables. Illustrons maintenant
cette d\'efinition (non formelle) par un exemple. Soient $A = \q X \{X\f
\q
Y (Y\f X) \}$, $\{X_i/i \geq 0 \}$ et $\{Y_i/i \geq 0 \}$ deux ensembles
disjoints de variables propositionnelles. Les sous-types n\'egatifs de
$A$ sont
: $X_i, Y_i, i \geq 0$ et ses sous-types positifs sont : $A$, $\q Y (Y\f
X_i), Y_j\f X_i, X_i$ avec $i, j \geq 0$.\\

D'autre part, on appelle {\bf type avec substitution}, tout type de la
forme
 $A [t_1/x_1,..., t_n/x_n]$, obtenu par substitution simultan\'ee de
  $t_1$ \`a $x_1$,..., $t_n$ \`a $x_n$ dans le type $A$, o\`u $x_1$,...,
$x_n$ sont des
  variables d'individu, et $t_1$,..., $t_n$ des termes du langage.\\

On dit qu'un type $\q^+$ (resp. $\q^-$) $A$ {\bf satisfait la condition
(*)} si,
lorsque $B$ et $C$ sont deux sous-types n\'egatifs (resp. positifs) avec
substitution de $A$,
alors on ne peut pas
avoir les
propri\'et\'es suivantes :
\begin{center}
 $B < G$, $G \sim C_n$, \\
et\\
$C \leq C_1 \f (C_2 \f (... \f (C_n \f D)...))$,
\end{center}
o\`u $C_1, C_2,..., C_n, D$ et $G$ sont des types du syst\`eme ${{\cal
AF}2}$.\\
Rappelons bien que $B < G$ signifie que $B \leq G$ et $B \neq G$.\\

Pour la v\'erification de la condition (*), il s'agit d'un test
s'appartenant \`a l'unification
modulo une th\'eorie \'equationnelle $E$ (ou $E$-unification), et donc
il est facile de voir
 qu'il est ind\'ecidable lui aussi,
 par exemple par recodage du probl\`eme de Post ou du dixi\`eme
probl\`eme
  de Hilbert (voir [3]). \\
   
 On dit qu'un type $A$ est un {\bf bon type positif} (resp. {\bf
bon type n\'egatif}),
not\'e en abr\'eg\'e ${\cal B^+}$ (resp. ${\cal B^-}$) s'il est $\q^+$
(resp.
$\q^-$) et
satisfait la condition (*).\\

{\bf Remarques.} (1) J.-L. K\textsc{rivine} a pr\'esent\'e dans [6] une
m\'ethode pour
d\'efinir les types
de donn\'ees
syntaxiques du syst\`eme ${{\cal AF}2}$. Par exemples :
\begin{itemize}
\item Le type bool\'een est la formule :
\begin{center}
 $B [x] = \q X\{X0 \f (X1 \f Xx)\}$,
\end{center}
 o\`u $0$ et $1$ sont des symboles de constante.
 \item Le type des entiers naturels est d\'ej\`a d\'efini.
\item  Le type des listes d'\'el\'ements de type $U$ est la formule :
\begin{center}
 $LU [x] = \q X\{X \emptyset \f [\q y \q z (U [y]\f (Xz\f Xcons(y,
z)))\f
Xx]\}$,
\end{center}
 o\`u $cons$ est un symbole de fonction binaire et $\emptyset$ une
constante.
\end{itemize}

On peut v\'erifier que ces types de donn\'ees sont des types ${\cal
B^+}$. \\ 

(2) K. N\textsc{our} a d\'efini dans [8] une autre classe des types de
donn\'ees
(types de donn\'ees descendants). On peut v\'erifier \'egalement que ces
types sont tous des ${\cal B^+}$.  \hfill $\spadesuit$ \\

Soient $L$ un langage du second ordre et $E$ un syst\`eme d'\'equations
de
$L$, ce qui d\'efinit le syst\`eme de typage ${{\cal AF}2}$.\\

On se propose de d\'emontrer le th\'eor\`eme suivant :\\

{\bf Th\'eor\`eme II Soient A un type ${\cal B^+}$ clos du syst\`eme
${{\cal AF}2}$, et $t$ un
$\l$-terme, alors:\\ $t \in \mid A \mid_{\b}$ $\so{ssi}$ $t \f_{\b} t'$
et
$\v_{{\cal AF}2} t' : A$.}\\

Pour d\'emontrer ce th\'eor\`eme, on a besoin de certaines d\'efinitions
et
d'un lemme.\\

Soit $A$ un type ${\cal B^+}$ du syst\`eme ${{\cal AF}2}$. On note
${\cal
B_s^+}$ (resp. ${\cal B_s^-}$), l'ensemble des sous-types positifs
(resp. n\'egatifs)
 avec substitution de $A$. Soient $\Omega = \{ x_i/i \in {\bf N} \}$ une
\'enum\'eration d'un ensemble
 infini de variables du $\l$-calcul, et notons par $A_i$ ($i \in {\bf
N})$
les sous-types ${\cal B_s^-}$ de $A$, o\`u chaque type $A_i$ se
r\'ep\`ete une
infinit\'e de fois. On d\'efinit alors l'ensemble $\mit
\G_A^- = \{x_i : A_i /i \in {\bf N} \}$. Soit $u$ un $\l$-terme, tel que
$Fv (u) \subseteq \Omega$. On d\'efinit le contexte $\mit \G_{A,u}^-$
comme
\'etant la restriction de $\mit \G_A^-$ sur les d\'eclarations contenant
les variables de $Fv (u)$. La notation $\mit \G_A^- \v_{{\cal AF}2} u :
B$
exprime que $\mit \G_{A,u}^- \v_{{\cal AF}2} u : B$.\\
On pose $\mit \G_A^- \v_{{\cal AF}2}^{\b} u : B$ ssi il existe un
$\l$-terme $u'$, tel que
$u \f_{\b} u'$ et $\mit \G_A^- \v_{{\cal AF}2} u' : B$.\\

On d\'efinit le $\L_{\b}$-mod\`ele ${\cal M}$ et l'interpr\'etation
${\cal
I}$ de la m\^eme fa\c con qu'au
paragraphe $3$ (en rempla\c cant $\mit \G^-$ par $\mit \G_A^-$ et la $\b
\eta$-\'equivalence par la $\b$-r\'eduction). Alors on a le lemme
suivant
:\\

{\bf Lemme 4.1} {\it (i) Si $S$ est un sous-type ${\cal B_s^+}$ de $A$,
et
$\t \in \mid
S\mid_{\cal {M,I}}$, alors $\mit \G_A^-\v_{{\cal AF}2}^{\b} \t : S$.\\
(ii) Si $S$ est un sous-type ${\cal B_s^-}$ de $A$, et $\mit
\G_A^-\v_{{\cal
AF}2}^{\b} \t : S$, alors $\t \in \mid S\mid_{\cal {M,I}}$}.\\

{\bf Preuve} : Par induction simultan\'ee sur les sous-types ${\cal
B_s^+}$
et ${\cal B_s^-}$ de $A$ (la complexit\'e \'etant le nombre de symboles
logiques dans le type).\\

\so{Preuve de (i)}
\begin {itemize}
\item Si $S$ est atomique, alors on a la m\^eme preuve que
celle du lemme 3.2.
\item Si $S = \q XB$, o\`u $B$ est ${\cal B_s^+}$ , alors $\t\in \mid \q
XB
\mid_{\cal {M,I}} \Leftrightarrow (\q \Phi : \mid {\cal M} \mid^n \f
P_{S_{\b}(\L)}) (\t\in \mid B[X] \mid_{{{\cal M,I} [X \leftarrow
\Phi]}}$.
 Soit $Y$ une variable de relation $n$-aire qui ne figure pas dans $\mit
\G_{A,\t}^-$ et $B$. Donc $\t\in \mid B[X] \mid_{{{\cal M,I} [X
\leftarrow \mid Y
\mid_{\cal {M,I}}]}} = \mid B[Y / X] \mid_{\cal {M,I}}$, d'apr\`es le
lemme
2.3. D'o\`u par hypoth\`ese d'induction, $\mit \G_A^- \v_{{\cal
AF}2}^{\b} \t
: B[Y]$, donc il existe un $\l$-terme $\t'$ tel que $\t \f_{\b} \t'$ et
$\mit \G_{A,\t'}^- \v_{{\cal AF}2} \t' :
B[Y]$. Comme $Fv (\t') \subseteq Fv (\t)$, alors par le choix de $Y$, on
d\'eduit que  $\mit \G_{A,\t'}^- \v_{{\cal AF}2} \t' : \q YB[Y] = \q
XB$, et donc
$\mit \G_A^- \v_{{\cal AF}2}^{\b} \t : S$.
\item Si $S = B\f C$, o\`u $B$ est ${\cal B_s^-}$ et $C$ est ${\cal
B_s^+}$, alors soit $\t \in \mid B\f C \mid_{\cal {M,I}}$, et soit $y$
une
variable du $\l$-calcul telle que $y : B \in \mit \G_A^-$. On a $y : B
\v_{{\cal AF}2} y : B$, donc, d'apr\`es (ii), $y \in \mid B \mid_{\cal
{M,I}}$, par suite $(\t)y \in \mid C \mid_{\cal {M,I}}$, et donc,
d'apr\`es
l'hypoth\`ese d'induction, $\mit \G_A^- \v_{{\cal AF}2}^{\b}
(\t)y : C$. D'o\`u il existe un $\l$-terme $\t'$ tel que $(\t)y \f_{\b}
\t'$, et $\mit \G_{A,\t'}^- \v_{{\cal
AF}2} \t': C$. Il en r\'esulte que $(\t)y$ est normalisable, et donc
$\t$
est normalisable. La forme normale de $\t$ est $x$ ou $(x)\t_1...\t_n$
$(n \geq 1)$ ou $\l x \th$.\\
{\so{\bf Cas 1}}: Si $\t\f_{\b} x$, alors $(\t)y\f_{\b} (x)y$.
Comme $\mit \G_{A,\t'}^-\v_{{\cal AF}2} \t' : C$, on aura $\mit
\G_{A,\t'}^- \v_{{\cal AF}2} (x)y : C$. Or $Fv ((x)y) \subseteq Fv
(\t')$,
donc $\mit \G_{A,(x)y}^- \v_{{\cal AF}2} (x)y : C$. D'o\`u, d'apr\`es
le corollaire 1.9, la variable $x$ est d\'eclar\'ee d'un type $F$ dans
$\mit \G_{A,(x)y}^-$, avec : $F \leq G\f D, D \sim  D', \mit
\G_{A,(x)y}^-
\v_{{\cal AF}2} y : G, C = \q \mbox{\boldmath$\xi$} D'$, et
\mbox{\boldmath$\xi$} n'est pas libre dans $\mit \G_{A,(x)y}^-$, o\`u
$G$ est
$\q^+$ et $D$ est $\q^-$ qui ne commence pas par des quantificateurs.\\
 $y : B \in \mit \G_{A,(x)y}^-$, donc $B \leq B', B' \sim B''$, $G = \q
\mbox{\boldmath$\xi'$} B''$, et
\mbox{\boldmath$\xi'$} n'est pas libre dans $\mit \G_{A,(x)y}^-$.\\
Supposons que \mbox{\boldmath$\xi$} commence par une variable du second
ordre $X$. Comme $X$ n'est pas libre dans $\mit \G_{A,(x)y}^-$, alors
$X$
n'est pas
libre dans $F$. Or $F \leq G \f D$, ce qui implique que $X$ n'est pas
libre
dans $D$, donc non plus dans $D'$. Ce qui contredit le fait que $C$ est
propre.
 De plus \mbox{\boldmath$\xi$} ne peut pas commencer par une variable du
premier ordre, car $C$ est $\q^+$, d'o\`u $C = D'$.\\
Maintenant, comme $G$ est $\q^+$, on d\'emontre de la m\^eme mani\`ere
que $G = B''$.\\
D'autre part on a $B = B'$, car sinon, i.e si $B < B'$, alors, comme $F
\leq G \f D$ et $G \sim
B'$, on obtient une contradiction avec la condition (*).\\
Finalement on a $\mit \G_{A,x}^- \v_{{\cal AF}2} x : F$, par
cons\'equent $\mit
\G_{A,x}^- \v_{{\cal AF}2} x :
G\f D$, et donc $\mit \G_{A,x}^- \v_{{\cal AF}2} x : G\f D'$. D'o\`u
$\mit
\G_{A,x}^-
\v_{{\cal AF}2} x : B''\f C$, donc $\mit \G_{A,x}^- \v_{{\cal AF}2} x :
B\f
C$, et par suite $\mit \G_A^- \v_{{\cal AF}2}^{\b} \t : S$.\\
{\so{\bf Cas 2}}: Si $\t\f_{\b} (x)\t_1...\t_n$, alors
$(\t)y\f_{\b}(x)\t_1...\t_ny$. Comme $\mit \G_{A,\t'}^- \v_{{\cal AF}2}
\t' :
C$ et $Fv ((x)\t_1...\t_ny) \subseteq Fv (\t')$, on  d\'eduit que
$\mit \G_{A,(x)\t_1...\t_ny}^- \v_{{\cal AF}2} (x)\t_1...\t_ny : C$.
Donc,
d'apr\`es le corollaire 1.9, $x$ est d\'eclar\'ee d'un type $F$ dans
$\mit
\G_{A,(x)\t_1...\t_ny}^-$ avec : $F \leq C_1\f D_1$, $D_1\sim D'_1$,
$D'_i
= C_{i+1}\f D_{i+1}$,
$D_{i+1}\sim D'_{i+1}$ ($1 \leq i \leq n-1$), $D'_n = G\f D_{n+1}$,
$D_{n+1}\sim
D'_{n+1}$,
$\mit \G_{A,(x)\t_1...\t_ny}^- \v_{{\cal AF}2}
\t_i : C_i$ ($1 \leq i \leq n$), $\mit \G_{A,(x)\t_1...\t_ny}^-
\v_{{\cal
AF}2} y : G$, $C = \q
\mbox{\boldmath$\xi$} D'_{n+1}$, et
\mbox{\boldmath$\xi$} n'est pas libre dans $\mit
\G_{A,(x)\t_1...\t_ny}^-$.\\
 $y : B \in \mit \G_{A,(x)\t_1...\t_ny}^-$, donc $B\leq B', B'\sim B''$,
$G = \q
\mbox{\boldmath$\xi'$} B''$, et \mbox{\boldmath$\xi'$} n'est pas libre
dans
$\mit \G_{A,(x)\t_1...\t_ny}^-$.\\
\mbox{\boldmath$\xi$} ne commence pas par une variable du premier ordre,
car $C$ est $\q^+$.\\
Si \mbox{\boldmath$\xi$} commence par une variable du second ordre $X$,
et
comme $C$ est propre, cette variable doit \^etre libre dans $D'_{n+1}$,
donc elle est libre dans $D_{n+1}$, et donc dans $D'_n$ aussi. Par
cons\'equent $X$ est libre dans tous les $D_i$ (resp. $D'_i$) ($1 \leq i
\leq n + 1$). Donc $X$ est libre dans $F$, et par suite dans $\mit
\G_{A,(x)\t_1...\t_ny}^-$, ce qui est impossible. D'o\`u $C =
D'_{n+1}$.\\
De la m\^eme mani\`ere et comme $G$ est $\q^+$, on d\'emontre que $G =
B''$.\\
D'autre part, on remarque que $B = B'$, car sinon, i.e si $B < B'$, et
comme on a montr\'e que $F \leq C_1 \f (C_2 \f (... \f (C_n \f (G \f
C)...)))$ et $B' \sim G$,
alors on obtient une contradiction, d'apr\`es la
condition (*).\\
Finalement on a $\mit \G_{A,(x)\t_1...\t_n}^- \v_{{\cal AF}2} x : F$,
donc
$\mit \G_{A,(x)\t_1...\t_n}^-
\v_{{\cal AF}2} (x)\t_1...\t_n :
G\f D_{n+1}$. Par cons\'equent $\mit \G_{A,(x)\t_1...\t_n}^- \v_{{\cal
AF}2} (x)\t_1...\t_n :
G\f D'_{n+1}$, et alors
$\mit \G_{A,(x)\t_1...\t_n}^- \v_{{\cal AF}2} (x)\t_1...\t_n : G\f C$.
D'o\`u $\mit \G_{A,(x)\t_1...\t_n}^-
\v_{{\cal AF}2} (x)\t_1...\t_n : B\f C$, et donc $\mit \G_A^- \v_{{\cal
AF}2}^{\b} \t : S$\\
 {\so{\bf Cas 3}}: Si $\t\f_{\b}\l x \th$, alors, comme l'ensemble $\mit
\G^-_A$ contient une infinit\'e de d\'eclarations pour chaque sous-type
 ${\cal B_s^-}$ de $A$, soit $y$ une variable d\'eclar\'ee de type $B$
dans
 $\mit \G^-_A$, n'appartenant pas \`a $Fv (\th)$. Alors $(\t)y
\f_{\b}(\l x \th)y \f_{\b} \th[y/x]$.
On a $\mit \G_{A,\t'}^- \v_{{\cal AF}2} \t' : C$, donc $\mit
\G_{A,\t'}^-
\v_{{\cal AF}2} \th[y/x] : C$, et donc, $\mit \G_{A,\th[y/x]}^-
\v_{{\cal
AF}2} \th[y/x] : C$, car $Fv (\th[y/x]) \subseteq Fv (\t')$. D'o\`u
$\mit
\G_{A,\l y \th[y/x]}^- \v_{{\cal AF}2} \l y \th[y/x] : B\f C$. Comme $y$
n'est pas libre
dans $\th$, alors $\l y \th[y/x] = \l x \th$, par cons\'equent
 $\mit \G_A^- \v_{{\cal AF}2} \l x \th : S$ et $\mit \G_A^- \v_{{\cal
AF}2}^{\b} \t : S$.
\end {itemize}

\so{Preuve de (ii)}
\begin{itemize}
\item Si $S$ est atomique, alors on reprend la m\^eme preuve de (ii) du
lemme 3.2.
\item Si $S = B \rightarrow C$, o\`u $B$ est ${\cal B_s^+}$  et $C$ est
${\cal B_s^-}$,
supposons $\mit \G_A^- \v_{{\cal AF}2}^{\b} \t : B \rightarrow C$.
 Donc il existe un $\l$-terme $\t'$ tel que $\t \f_{\b} \t'$ et $\mit
\G_{A,\t'}^- \v_{{\cal AF}2} \t' : B \rightarrow C$.
Si $u \in \mid B \mid_{\cal {M,I}}$,
alors d'apr\`es (i), $\mit \G_A^-\v_{{\cal AF}2}^{\b}  u : B$, donc il
existe un
 $\l$-terme $u'$ tel que $u
\f_{\b} u'$ et $\mit \G_{A,u'}^- \v_{{\cal AF}2} u' : B$. D'o\`u $\mit
\G_{A,(\t')u'}^-\v_{{\cal AF}2} (\t') u' : C$, et comme $(\t)u \f_{\b}
(\t')u'$, alors $\mit \G_A^-\v_{{\cal AF}2}^{\b} (\t)u : C$. D'o\`u,
d'apr\`es l'hypoth\`ese d'induction, $(\t)u \in \mid C \mid_{\cal
{M,I}}$.
Par cons\'equent $\t \in \mid B \f C \mid_{\cal {M,I}}$.
\item Si $S = \q xB$, o\`u $B$ est ${\cal B_s^-}$, supposons que $\mit
\G^-_A \v_{{\cal AF}2}^{\b} \t : \q xB$, et soit $a \in \mid {\cal M}
\mid$. Alors on a : $a = \sou{b}$, o\`u $b$ est un terme de $L$, et il
existe 
un $\l$-terme $\t'$ tel que $\t \f_{\b} \t'$ et $\mit
\G_{A,\t'}^-\v_{{\cal AF}2} \t' : \q xB$, par suite
$\mit \G_{A,\t'}^-\v_{{\cal AF}2} \t' : B [b/x]$. Or, par d\'efinition,
$B[b/x]$ est un sous-type ${\cal B_s^-}$ de $A$, donc d'apr\`es
l'hypoth\`ese d'induction, $\t' \in \mid B[b/x] \mid_{\cal {M,I}} = \mid
B
\mid_{{{\cal M,I} [x \leftarrow \sou{b}]}} = \mid B \mid_{{{\cal M,I} [x
\leftarrow a]}}$ (d'apr\`es le lemme 2.2). Par cons\'equent $\t \in \mid
B
\mid_{{{\cal M,I} [x \leftarrow a]}}$, et ce pour tout $a \in \cal M$;
donc
 $\t \in \mid \q xB \mid_{\cal M,I}$. \hfill $\spadesuit$
\end{itemize}

On peut d\'eduire ainsi la preuve du th\'eor\`eme II.\\

{\bf Preuve du th\'eor\`eme II }
\begin{itemize}
\item La condition suffisante est une cons\'equence directe du
th\'eor\`eme 2.1,
et du fait que l'interpr\'etation d'un type du syst\`eme ${{\cal AF}2}$
dans un
$\L_{\b}$-mod\`ele $M$ est une partie ${\f_\b}$-satur\'ee.
\item La condition est n\'ecessaire. En effet : Soient $A$ un type
${\cal
B^+}$ clos, et
$t$ un $\l$-terme tel que $t \in \mid A \mid_{\b}$, alors $t \in \mid
A \mid_{\cal M}$, pour tout $\L_{\b}$-mod\`ele ${\cal M}$ associ\'e \`a
un
ensemble $\mit \G_A^-$ (comme d\'ecrit avant). De plus on peut supposer
que
$\mit \G_A^-$ ne contient pas de d\'eclarations pour les variables
libres
de $t$, i.e $\mit \G_{A,t}^- = \emptyset$. D'apr\`es le (i) du lemme
4.1,
$\mit \G_A^-\v_{{\cal AF}2}^{\b} t : A$, donc il existe un $\l$-terme
$t'$ tel
que $t \f_{\b} t'$ et $\mit
\G_{A,t'}^-\v_{{\cal AF}2} t' : A$. Comme $Fv (t') \subseteq Fv (t)$,
alors
$\mit \G_{A,t'}^- = \emptyset$, d'o\`u le r\'esultat.    \hfill
$\spadesuit$\\
\end {itemize}

 D'apr\`es le th\'eor\`eme 3.1, on peut d\'eduire le th\'eor\`eme
suivant :\\

{\bf Th\'eor\`eme 4.2} {\it Soient A un type ${\cal B^+}$ clos du
syst\`eme
${{\cal AF}2}$, et $t$ un
$\l$-terme, alors :\\ $t \in \mid A \mid_{\b}$ $\so{ssi}$ il existe un
$\l$-terme
 $t'$ tel que $t \f_{\b} t'$ et $\v_{{\cal AF}2_0} t' : A$.}\\

{\bf Corollaire 4.3} {\it Soient $A$ un type ${\cal B^+}$ clos du
syst\`eme
${{\cal AF}2}$, et $t$ un
$\l$-terme. \\(i) Si $t \in \mid A \mid_{\b}$,
alors $t$ est normalisable et se r\'eduit par $\b$-r\'eduction \`a un
terme clos.\\
(ii) $\mid A \mid_{\b}$ est stable par $\b$-\'equivalence (i.e si $t \in
\mid A \mid_{\b}$ et $t \simeq\sb{\b} t'$, alors $t' \in \mid A
\mid_{\b}$).}\\

{\bf Preuve} : (i) Si $t \in \mid A \mid_{\b}$, alors,
d'apr\`es le th\'eor\`eme II, il existe un $\l$-terme $t'$ tel que $t
\f_{\b} t'$ et $\v_{{\cal AF}2} t' :
A$. Donc $t$ se r\'eduit par $\b$-r\'eduction \`a un terme normalisable
clos, d'o\`u le r\'esultat.\\
(ii) Soit $t \in \mid A \mid_{\b}$, avec $t \simeq\sb{\b} t'$, alors il
existe un $\l$-terme $v$ tel que $t\f_{\b} v$ et $t'\f_{\b} v$.
Or d'apr\`es le th\'eor\`eme II, il existe un $\l$-terme $u$ tel que $t
\f_{\b} u$ et $\v_{{\cal AF}2} u : A$. On peut supposer que $u$  est
normal
(car sinon $u\f_{\b} w$, $w$ normal et $\v_{{\cal AF}2} w : A$). D'o\`u
$v
\f_{\b} u$, et comme $\v_{{\cal AF}2} u : A$, alors d'apr\`es le
th\'eor\`eme 2.1, $u \in \mid A \mid_{\b}$. Mais $\mid A \mid_{\b}$ est
$\f_{\b}$-satur\'ee, il en r\'esulte que $v \in \mid A \mid_{\b}$ et par
suite
$t' \in \mid A \mid_{\b}$. \hfill $\spadesuit$\\

{\bf Remarques.} Nous allons voir que les conditions qui d\'efinissent
un
type ${\cal B^+}$ sont toutes n\'ecessaires pour avoir le th\'eor\`eme
II.
\\

 (1) Consid\'erons les types
\begin{eqnarray}
A &= &\q X \{ (X0\f \q yXy)\f (X0\f X0) \} \nonumber \\
B &= &\q X \{ \q x (X0 \f (Xx \f X0)) \f (X0 \f \q x (Xx \f X0)) \}
\nonumber  \\
C &= &\q X \{ ( \q xXx \f X0) \f (X0\f X0) \} \nonumber
\end{eqnarray}
Il est clair que ces types sont propres, satisfont la condition (*),
mais ils ne sont pas $\q^+$ (dans $A$ et $B$ on a un quantificateur du
premier ordre derri\`ere
une fl\`eche, et dans $C$, le probl\`eme vient d'un quantificateur du
premier ordre positif). \\

 $\not \v_{{\cal AF}2} I = \l xx : A$, car sinon on
aura, $x : X0 \f \q yXy \v_{{\cal AF}2} x : X0 \f X0$. Ce qui est
impossible, d'apr\`es le
th\'eor\`eme 1.8.\\ Pourtant $I \in \mid A \mid_{\b}$ : Soit $M$ un
$\L_{\b}$-mod\`ele quelconque, il faut montrer que $I \in \mid (X0
 \f \q yXy) \f (X0 \f X0)\mid_{M, [X \leftarrow \Xi]}$ pour tout $\Xi
\in
P_{S_{\b}(\L)}^{\mid M \mid}$. Consid\'erons donc un \'el\'ement $\Xi$
dans
$P_{S_{\b}(\L)}^{\mid M \mid}$; si $u \in \mid(X0 \f \q yXy)\mid_{M, [X
\leftarrow \Xi]}$, et $v
\in \mid X0 \mid_M$. Alors on a : $(u)v \in \mid \q yXy\mid_{M, [X
\leftarrow \Xi]}$, et comme
$(I)uv \f_{\b} (u)v$, on voit
que $(I)uv \in \mid \q yXy\mid_{M, [X \leftarrow \Xi]}$,
 car $\mid\q yXy\mid_{M, [X \leftarrow \Xi]}$ est une
partie ${\f_\b}$-satur\'ee de $\L$. Par cons\'equent
 $(I)uv \in \mid X (0)\mid_{M, [X \leftarrow \Xi]} =
\Xi (0_M)$.\\

De m\^eme pour $B$, il est clair que $\not \v_{{\cal AF}2} I : B$. De
plus, si $M$ est un
$\L_{\b}$-mod\`ele et $\Xi \in P_{S_{\b}(\L)}^{\mid M \mid}$, alors en
supposant que $u \in \mid \q x (X0 \f (Xx \f X0)) \mid_{M, [X \leftarrow
\Xi]}$,
 $v \in \mid X0 \mid_{M, [X \leftarrow \Xi]}$, et en prenant un
\'el\'ement quelconque
$a$ de $\mid M \mid$, on voit que $u \in \mid (X0 \f (Xa \f X0))
\mid_{M, [X \leftarrow \Xi]}$,
et donc $(u)v \in \mid Xa \f X0 \mid_{M, [X \leftarrow \Xi]}$,
 d'o\`u $(I)uv \in \mid Xa \f X0 \mid_{M, [X \leftarrow \Xi]}$. Par
cons\'equent $I
\in \mid B \mid_M$.\\

On reprend la m\^eme d\'emonstration pour $C$.\\

(2) Soit
\begin{center}
$C' = \q X \{ (\q YX0 \f X0) \f (X0 \f X0) \}$
\end{center}
$C'$ satisfait la condition (*), et v\'erifie toutes les propri\'et\'es
qui
figurent dans la
d\'efinition d'un type $\q^+$, sauf le fait qu'il n'est pas propre. On
peut
d\'emontrer comme dans
(1), que $I \in \mid C' \mid_{\b}$, et $\not \v_{{\cal AF}2} I : C'$.\\

(3) Reprenons l'exemple du type $D = \q X\{\q Y (Y \f X) \f X\}$
d\'ecrit
dans le paragraphe 3. $D$
est un type propre qui satisfait la condition (*), mais il n'est pas ni
$\q^+$ ni $\q_2^+$. Par
une d\'emonstration analogue \`a celle qui est d\'ej\`a faite, on voit
que le
$\l$-terme $t = \l x (x) (\delta)\delta \in \mid D \mid_{\b}$, et $t$
n'est pas normalisable.\\
De m\^eme $t'=\l x (x) y \in \mid D \mid_{\b}$, et $t'$ est normal et
non clos.\\

(4) Consid\'erons enfin les types
\begin{eqnarray}
E &= &\q X \{ \q x (Xx \f X0) \f (\q xXx \f X0) \} \nonumber \\
F &= &\q X \{ \q x (X0 \f (Xx \f X0)) \f (X0 \f (\q xXx \f X0)) \}
\nonumber \\
K &= &\q X\q Y \{\q y\{[\q x(X(x,0)\f X(0,0))\f (\q xX(x,y)\f X(0,0))]\f
Y\}\f Y\} \nonumber.
\end{eqnarray}
$E$ et $F$ sont $\q^+$, mais ils ne satisfont pas la condition (*). En
effet : \\
Pour $E$, on a : $\q x (Xx \f X0) \leq (Xx \f X0)$, et $\q xXx < Xx$. \\
Et pour $F$, on a : $\q x (X0 \f (Xx \f X0)) \leq X0 \f (Xx \f X0)$, et
$\q
xXx < Xx$.\\
 De la m\^eme mani\`ere, on d\'emontre que $I \in \mid E \mid_{\b}$
(resp. $I
\in \mid F \mid_{\b}$), et que $\not \v_{{\cal AF}2} I : E$ (resp. $\not
\v_{{\cal AF}2} I : F$).\\

 $K$ n'est pas ${\cal B^+}$ pour des raisons plus compliqu\'ees. En
effet :
$\q x(X(x,0)\f X(0,0)) \leq X(x,0)\f X(0,0)$ et $\q xX(x,y) [0/y] <
X(x,0)$.\\
On d\'emontre facilement que $\l x (x)I \in \mid K \mid_{\b}$ et $\not
\v_{{\cal AF}2} \l x (x)I : K$.  \hfill $\spadesuit$\\

Il est clair qu'un type du syst\`eme ${\cal F}$ est ${\cal B^+}$ ssi il
est $\q_2^+$ et propre. D'o\`u le r\'esultat suivant:\\

{\bf Th\'eor\`eme 4.4} {\it Soient A un type $\q_2^+$, clos et propre du
syst\`eme ${\cal F}$, et $t$ un
$\l$-terme, alors:\\  $t \in \mid A \mid_{\b}$  $\so{ssi}$ il existe un
$\l$-terme
$t'$ tel que $t
\f_{\b} t'$ et $\v_{\cal F} t' : A$.}\\

{\bf Remerciements}. Nous remercions C. R\textsc{affalli} pour ses
remarques.\\

\end{document}